\documentclass[reqno,11pt]{article}
\usepackage{amsmath,amssymb,epsfig}
\usepackage{mathrsfs,enumerate}
\usepackage{amsbsy,bm}

\usepackage{theorem}
\newtheorem{theorem}{Theorem}[section]  
\newtheorem{lemma}[theorem]{Lemma}
\newtheorem{corollary}[theorem]{Corollary}
\newtheorem{proposition}[theorem]{Proposition}

\theorembodyfont{\rmfamily}
\newtheorem{remark}[theorem]{Remark}
\newtheorem{example}[theorem]{Example}
\newtheorem{definition}[theorem]{Definition}

\DeclareMathAlphabet{\mathpzc}{OT1}{pzc}{m}{it}

\newcommand{\prob}[1]{\mathscr P(#1)}
\newcommand{\probt}[1]{\mathscr P_2(#1)}

\newcommand{\la}{\left\langle}
\newcommand{\ra}{\right\rangle}
\newcommand{\lims}{\varlimsup}
\newcommand{\limi}{\varliminf}
\newcommand{\restr}[1]{\lower3pt\hbox{$|_{#1}$}}

\newenvironment{Proof}{\removelastskip\par\medskip   
\noindent{\em Proof.} \rm}{\penalty-20\null\hfill$\square$\par\medbreak}

\newenvironment{Proofof}{\removelastskip\par\medskip   
\noindent{\em Proof of} \rm}{\penalty-20\null\hfill$\square$\par\medbreak}

\makeatletter \@addtoreset{equation}{section} \makeatother

\newcommand{\e}{\mathrm{e}}
\newcommand{\R}{\mathbb{R}}
\newcommand{\N}{\mathbb{N}}
\newcommand{\ep}{\varepsilon}

\newcommand{\abra}[1]{\left( #1 \right)}

\newcommand{\gl}[1]{| \nabla\!_L #1 |}
\newcommand{\ls}[1]{| \nabla\!_- #1 |}

\newcommand{\gs}[1]{| \nabla #1 |}

\newcommand{\gss}[2]{\la \nabla #1 , \nabla #2 \ra}

\newcommand{\lap}{\Delta}

\newcommand{\vol}{\mathcal{H}^n}

\newcommand{\pacs}[1]{\mathscr P^*(#1)}

\newcommand{\bbM}{\mathbb{M}}
\newcommand{\CD}{\mathsf{CD}}

\newcommand{\Ent}{\mathop{\rm Ent}\nolimits}
\newcommand{\Ric}{\mathop{\rm Ric}\nolimits}

\title{Heat flow on Alexandrov spaces}
\author{Nicola Gigli\footnote{Institut f\"ur Angewandte Mathematik,
Universit\"at Bonn, Endenicher Allee 60, 53115 Bonn, Germany 
/ current affiliation: Department of Mathematics, University of Nice, 
Valrose 06108 Nice Cedex 02, France
(\texttt{gigli@unice.fr})}
\footnote{Partially supported by SFB 611},
Kazumasa Kuwada\footnote{Graduate School of Humanities and Sciences,
Ochanomizu University, Ohtsuka 2-1-1, Bunkyo-ku, Tokyo 112-8610, Japan
\& Institut f\"ur Angewandte Mathematik, Universit\"at Bonn,
Endenicher Allee 60, 53115 Bonn, Germany
(\texttt{kuwada.kazumasa@ocha.ac.jp})}
\footnote{Partially supported by the JSPS fellowship for research abroad
\& the Grant in Aid for Young Scientists (B) 22740083},
Shin-ichi Ohta\footnote{Department of Mathematics, Kyoto University,
Kyoto 606-8502, Japan \& Max-Planck-Institut f\"ur Mathematik, Vivatsgasse 7,
53111 Bonn, Germany  (\texttt{sohta@math.kyoto-u.ac.jp})}
\footnote{Partially supported by the Grant-in-Aid for Young Scientists (B) 20740036}}
\date{\today}

\begin{document}

\maketitle

\begin{abstract}
We prove that on compact Alexandrov spaces with curvature bounded below
the gradient flow of the Dirichlet energy in the $L^2$-space produces
the same evolution as the gradient flow of the relative entropy
in the $L^2$-Wasserstein space. 
This means that the heat flow is well defined by either one of the two gradient flows. 
Combining properties of these flows, 
we are able to deduce the Lipschitz continuity of the heat kernel 
as well as Bakry-\'Emery gradient estimates and the $\Gamma_2$-condition. 
Our identification is established by purely metric means, 
unlike preceding results relying on PDE techniques. 
Our approach generalizes to the case of heat flow with drift.
\end{abstract}

\begin{section}{Introduction}

The heat equation is one of the most important evolutionary PDEs.
It is a well known fact in modern analysis that such an equation,
say in $\R^n$, can be seen as the gradient flow of the Dirichlet energy
\begin{equation}\label{eq:dir}
\frac{1}{2}\int_{\R^n} |\nabla f|^2 \,dx
\end{equation}
in the space $L^2(\R^n,\mathcal L^n)$.
This viewpoint has been extended to the concept of Dirichlet form 
(see e.g. \cite{BH,FOT}) and it has grown up into a huge research field 
in potential analysis and probability theory. 
More recently, Jordan, Kinderlehrer and Otto~\cite{JKO:FP:98}
understood that the same equation can be seen as the gradient flow
of the relative entropy 
\begin{equation}\label{eq:ent}
\int_{\R^n} \rho\log\rho \,dx
\end{equation}
in the $L^2$-Wasserstein space $(\probt{\R^n},W_2)$,
where $\probt{\R^n}$ stands for the space of Borel probability
measures on $\R^n$ with finite second moment.
This intuition, with the further studies of Otto~\cite{otto:geom:01},
has been one of the fundamental ingredients that drove the research
in the field of gradient flows in relation with optimal transportation
problems in the past decade (see, e.g., \cite{AGS} and \cite{Villani08}).

The aim of this paper is to carry on the study of the heat flow as gradient flow
of the two very different functionals, \eqref{eq:dir} and \eqref{eq:ent}
in the two metric spaces $L^2(\R^n,\mathcal L^n)$ and $(\probt{\R^n},W_2)$,
in a non-smooth setting.
The point is the following.
On the one hand, it is known that these two gradient flows produce the same evolution
`in all the smooth settings', i.e., it has been proved that in $\R^n$,
in Riemannian manifolds with Ricci curvature bounded below
and in compact Finsler manifolds (for the latter see the work of Sturm and
the third author~\cite{heatfinsl}), the gradient flow of the Dirichlet energy
with respect to the $L^2$-distance coincides with the gradient flow of
the relative entropy with respect to the Wasserstein distance $W_2$
(see also \cite{FSS}, \cite{Jui} for related work on different kinds of spaces).
On the other hand, these flows are studied well also in some
non-smooth settings.
These raise the natural question: do these two notions always coincide?

We remark that a natural abstract setting where one could try to
give an answer to this question is the one of metric measure spaces
with Ricci curvature bounded below (defined by Sturm~\cite{SturmI}
and Lott and Villani~\cite{LV-CD}, see Definition~\ref{df:CD}).
Indeed, in order to define either one of the two gradient flows,
one usually needs both a metric and a measure on the space considered.
Furthermore, it is known that in the Riemannian setting the fact that
the heat flow does not lose mass is strictly related to the bound
from below on the Ricci curvature, so that one has to assume
some kind of lower Ricci curvature bound also in the non-smooth setting.

In this paper, we restrict our attention to the case of finite dimensional,
compact Alexandrov spaces $(X,d,\vol)$ of curvature bounded below
without boundary,
equipped with the Hausdorff measure $\vol$.
An Alexandrov space of curvature bounded from below by $k$ with $k \in \R$
is a metric space of sectional curvature bounded from below by $k$ in the sense of
the triangle comparison property (see Definition~\ref{df:Alex}).
As naturally expected and recently shown in \cite{PetCD} and \cite{ZZ1},
such spaces actually have Ricci curvature bounded below
in the sense of Lott-Sturm-Villani (see Remark~\ref{rem:CD}).
Our following main result is the first one establishing the equivalence
of the two gradient flows in a genuinely non-smooth setting
(see Theorem~\ref{thm:main} for the slightly more general statement):

\begin{theorem}\label{thm:intro}
Let $(X,d,\vol)$ be a finite-dimensional compact Alexandrov space without boundary.
Then for any $f_0\in L^2(X,\vol)$ with $f_0\vol \in \prob{X}$,
the gradient flow $(f_t)_{t \in [0,\infty)}$ of the Dirichlet energy
with respect to the $L^2$-distance gives
the gradient flow $(f_t \vol)_{t \in [0,\infty)}$
of the relative entropy with respect to the $L^2$-Wasserstein distance $W_2$,
and vice versa.
\end{theorem}
On Alexandrov spaces, the structure of the heat flow as gradient flow
has already been studied in \cite{Kuw_Mac_Shi} for the Dirichlet energy approach
and in \cite{Ohta_grad-Alex} by the third author (see also \cite{Gigli-Ohta10} and \cite{Sav})
for the relative entropy approach.
However, up to now it was not clear that these two notions coincide.

We observe that the idea behind our proof of the theorem is
completely different from the ones used in the aforementioned smooth settings.
Indeed, all the proofs available in the smooth setting had the following structure:
\begin{itemize}
\item one studies the gradient flow of the Dirichlet energy with respect to $L^2$
and writes down the equation that it satisfies;
\item one studies the gradient flow of the relative entropy with respect to $W_2$
and writes down the equation that it satisfies;
\item one realizes that the two equations are actually the same and calls into
play the uniqueness of the solution of the equation to conclude that the two gradient flows coincide.
\end{itemize}
Our approach, instead, is completely different and in particular does not
pass from the study of the heat equation in the Alexandrov setting.
Read back in $\R^n$ our proof gives a new and purely metric way
to prove the coincidence of these two gradient flows. 

After having proved the identification, we can combine them: 
their interaction gives fruitful applications.
Among them, the Lipschitz continuity (Theorem~\ref{thm:HK-Lip}) of the heat kernel 
is derived as a corollary of a Bakry-\'Emery type gradient estimate 
(Theorem~\ref{thm:BE}), 
which follows with the aid of a result of the second author~\cite{K9}. 
The Lipschitz regularity improves the H\"older regularity
established in \cite{Kuw_Mac_Shi}. 
We believe that it is curious and worth to be underlined that such Lipschitz property 
immediately follows when one knows the equivalence of the two gradient flows, 
but is not at all trivial if one sticks to either one of the two approaches alone. 
Our Lipschitz regularity enables us to deduce the $\Gamma_2$-condition 
(Theorem~\ref{thm:Boch}) from a Lott-Sturm-Villani Ricci curvature bound.  
The $\Gamma_2$-condition as well as the Wasserstein contraction of heat flow 
and the Bakry-\'Emery gradient estimate for heat semigroup are known to be 
analytic characterizations of a lower Ricci curvature bound 
(see \cite{vRS,Bak97,Led_geom-Markov}). 
We show the equivalence of 
those three analytic conditions
even on Alexandrov spaces 
(with sharp constants, Theorem~\ref{thm:Ricci}). 

It should be remarked that 
our approach easily generalizes to the case of heat flow with drift,
where the drift is given by the gradient of a semiconvex potential $V$.
In other words, the current approach can be used to study the heat flow
on weighted Alexandrov spaces of the kind $(X,d,\e^{-V}\vol)$
(see the end of Section~\ref{sc:appl}).

After having completed the work on this paper, we got aware of a
paper~\cite{ZZ2} by Zhang and Zhu where the Lipschitz continuity
of the heat kernel has been studied, via a completely different argument.
In particular, their proof relies on the Lipschitz continuity of harmonic functions
due to Petrunin~\cite{PetERA}.
They even proceeded their study to \cite{ZZ3} and \cite{QZZ},
where they showed a Bochner type formula and Li-Yau estimates
under their own notion of lower Ricci curvature bound which is stronger
than the one of Lott-Sturm-Villani.
It is not discussed in \cite{ZZ2,ZZ3,QZZ} whether their approach generalizes 
to the case of weighted Alexandrov spaces. 

The organization of the article is as follows.
Section~\ref{sc:prel} is devoted to preliminaries for Alexandrov geometry
and known results on the gradient flow of the Dirichlet energy as well as
the gradient flow of the relative entropy on compact Alexandrov spaces.
We prove our main theorem in Section~\ref{sc:main},
and discuss its applications in Section~\ref{sc:appl}.
\end{section}

\begin{section}{Preliminaries and notations}\label{sc:prel}

\begin{subsection}{Alexandrov spaces}

We first review the basics of Alexandrov geometry,
see \cite{BGP}, \cite{Otsu-Shioya_JDG94} and \cite{BBI} for details.

Let $(X,d)$ be a metric space.
A rectifiable curve $\gamma:[0,l] \to X$ is called a geodesic if it is locally
length minimizing and parametrized with constant speed.
(Precisely, for any $t_0 \in [0,l]$, there is $\ep>0$ such that
$d(\gamma(s),\gamma(t))/|s-t|$ is constant for all
$s,t \in [t_0-\ep,t_0+\ep] \cap [0,l]$.)
If $\gamma$ is minimizing between endpoints, then we call it a minimal geodesic.
We say that $(X,d)$ is a geodesic space if any pair of points in $X$
are connected by a minimal geodesic.

For $k \in \R$, we denote by $\bbM^2(k)$ the simply-connected,
two-dimensional space form of constant sectional curvature $k$.
Given three points $x,y,z \in X$, with $d(x,y)+d(y,z)+d(z,x)<2\pi/\sqrt{k}$
if $k>0$, we can take a comparison triangle
$\triangle \tilde{x}\tilde{y}\tilde{z} \subset \bbM^2(k)$ such that
$d(\tilde{x},\tilde{y})=d(x,y)$, $d(\tilde{y},\tilde{z})=d(y,z)$ and
$d(\tilde{z},\tilde{x})=d(z,x)$.
Such a triangle is unique up to a difference of isometry.

\begin{definition}[Alexandrov spaces]\label{df:Alex}
For $k \in \R$, a complete geodesic space $(X,d)$ is called
an \emph{Alexandrov space of curvature bounded from below by $k$}
if, for any three points $x,y,z \in X$
(with $d(x,y)+d(y,z)+d(z,x)<2\pi/\sqrt{k}$ if $k>0$) and any minimal
geodesic $\gamma:[0,1] \to X$ from $y$ to $z$, we have
$d(x,\gamma(t)) \ge d(\tilde{x},\tilde{\gamma}(t))$ for all $t \in [0,1]$,
where $\triangle \tilde{x}\tilde{y}\tilde{z} \subset \bbM^2(k)$
is a comparison triangle of $\triangle xyz$ and
$\tilde{\gamma}:[0,1] \to \bbM^2(k)$ is the unique minimal geodesic
from $\tilde{y}$ to $\tilde{z}$.
\end{definition}

\begin{example}\label{ex:Alex}
(a) A complete Riemannian manifold is an Alexandrov space of curvature bounded from below by $k$
if and only if its sectional curvature is greater than or equal to $k$ everywhere.

(b) If $(X,d)$ is an Alexandrov space of curvature bounded from below by $k$,
then the scaled metric space $(X,c \cdot d)$ with $c>0$ is
an Alexandrov space of curvature bounded from below by $k/c^2$.

(c) For a convex domain $D$ in the Euclidean space $\R^n$,
the boundary $\partial D$ equipped with the length distance is
an Alexandrov space of nonnegative curvature.

(d) Let $(M,g)$ be a Riemannian manifold of nonnegative sectional curvature
and $G$ be a compact group acting on $M$ by isometries.
Then the quotient space $M/G$ equipped with the quotient metric
is an Alexandrov space of nonnegative curvature.

(e) If a sequence of Alexandrov spaces of curvature bounded from below by $k$ is convergent
with respect to the Gromov-Hausdorff distance,
then its limit space is again an Alexandrov space of curvature bounded from below by $k$.
\end{example}

Fix $x \in X$ and let $\hat{\Sigma}_x$ be the set of all unit speed geodesics
$\gamma:[0,l] \to X$ with $\gamma(0)=x$.
For $\gamma,\eta \in \hat{\Sigma}_x$, thanks to the curvature bound,
the joint limit
\[ \angle_x(\gamma,\eta) := \arccos \bigg( \lim_{s,t \to 0}
 \frac{s^2+t^2-d(\gamma(s),\eta(t))^2}{2st} \bigg) \]
exists and is 
a distance on $\hat{\Sigma}_x /\!\!\sim$ 
where $\gamma \sim \eta$ holds if $\angle_x(\gamma,\eta)=0$.
We define the space of directions $(\Sigma_x,\angle_x)$ at $x$
as the completion of $\hat{\Sigma}_x/\!\!\sim$ with respect to $\angle_x$. 

The tangent cone $(K_x,d)$ is the Euclidean cone
over $(\Sigma_x,\angle_x)$, i.e.,
\begin{align*}
K_x &:= \Sigma_x \times [0,\infty) / \Sigma_x \times \{0\}, \\
d \big( (\gamma,s),(\eta,t) \big)
&:= \sqrt{s^2+t^2-2st\cos\angle_x(\gamma,\eta)}.
\end{align*}
The inner product on $K_x$ is defined by
$\langle (\gamma,s),(\eta,t) \rangle_x :=st\cos\angle_x(\gamma,\eta)$.
In Riemannian manifolds, spaces of directions and tangent cones correspond to
unit tangent spheres and tangent spaces, respectively.

The Hausdorff dimension of $X$ is an integer or infinity.
From here on, \emph{we consider a compact $n$-dimensional
Alexandrov space of curvature bounded from below by $k$ without boundary
equipped with the $n$-dimensional Hausdorff measure $\vol$}
(see \cite{BGP} for the definition of the boundary of Alexandrov spaces).
We say that $x \in X$ is a singular point if $K_x$ is not isometric to the Euclidean space $\R^n$,
and denote the set of singular points by $S_X$.
We remark that $\vol(S_X)=0$ holds whereas $S_X$ can be dense in $X$
(see \cite{BGP}, \cite{Otsu-Shioya_JDG94}).
\end{subsection}

\begin{subsection}{Dirichlet energy and the associated gradient flow}

We introduce the Sobolev space and the Dirichlet energy following \cite{Kuw_Mac_Shi},
and will see that it coincides with other notions of Sobolev spaces. 
We begin by discussing a $C^1$-differentiable structure of the set of regular points
$X \setminus S_X$ established in \cite{BGP} and \cite{Otsu-Shioya_JDG94}.
We remark that Perelman extends this to $DC^1$-structure
(via `difference of concave functions', see \cite{Perel}),
but the $C^1$-structure is enough for considering the Sobolev space.
There is a weak $C^1$-atlas
$\{(U_{\phi},V_{\phi},\phi)\}_{\phi \in \Phi}$
in the sense that $U_{\phi} \subset X$ is an open set,
$\phi:U_{\phi} \to \R^n$ is a bi-Lipschitz embedding,
$V_{\phi} \subset U_{\phi}$ with
$\bigcup_{\phi \in \Phi} V_{\phi} \supset X \setminus S_X$,
and that the coordinate change $\phi_2 \circ \phi_1^{-1}$ is $C^1$ on
$\phi_1(V_{\phi_1} \cap V_{\phi_2} \cap (X \setminus S_X))$
if $V_{\phi_1} \cap V_{\phi_2} \neq \emptyset$
(\cite[Theorem~4.2(1)]{Otsu-Shioya_JDG94}).
Such charts are constructed through the distance function.
Precisely, $\phi$ is introduced as
$\phi(x):=(d(x,p_1),d(x,p_2),\ldots,d(x,p_n))$
for suitable $p_1,p_2,\ldots,p_n \in X$, and then $V_{\phi}$ is chosen
as the set of regular points $x$ such that a minimal geodesic
between $x$ and $p_i$ is unique for all $i$ (\cite[Section~3]{Otsu-Shioya_JDG94}).
It is worth mentioning that, for any $\ep>0$, $\phi$ can be
$(1+\ep)$-bi-Lipschitz by taking smaller $U_{\phi}$
(cf.\ \cite[Theorem~10.9.16]{BBI}).
We also remark that it is possible to modify $\phi$ by taking an average
so as to satisfy $V_{\phi}=U_{\phi}$ (\cite[Section~5]{Otsu-Shioya_JDG94}),
but it is unnecessary for our discussion (just like \cite[Remark~2.9]{Bert}).

We say that a function $f$ on $X$ is differentiable at a regular
point $x \in X \setminus S_X$ if $x \in V_{\phi}$ and $f \circ \phi^{-1}$
is differentiable at $\phi(x)$ for some $\phi \in \Phi$.
Then we can define the gradient vector $\nabla f(x) \in K_x$
by identifying $K_x$ and $\R^n$ through $\phi$.
(To be precise, by virtue of the first variation formula
(\cite[Theorem~3.5]{Otsu-Shioya_JDG94}),
each $d_i:=d(\cdot,p_i)$ is differentiable at $x$
with $\nabla d_i(x)=-v_i$, where $v_i \in K_x$ is the tangent vector
of the unique, minimal, unit speed geodesic from $x$ to $p_i$.
Then $K_x$ is linearly identified with $\R^n$ as $\sum_i a_i v_i=(-a_i)$.)
Moreover, again due to the first variation formula,
we obtain the Taylor expansion
\begin{equation}\label{eq:Taylor}
f\big( \gamma(t) \big)
 =f(x)+t\langle \nabla f(x),\dot{\gamma}(0) \rangle_x +o_x(t),
\end{equation}
where $\gamma:[0,\delta] \to X$ is a minimal geodesic emanating from $x$
and $o_x(t)$ is independent of the choice of $\gamma$
(see \cite[Lemma~3.4]{Bert}, and note that the remainder term
in the first variation formula for $d(\cdot,p_i)$ at $x$ indeed
depends only on $x$ and $p_i$).
Another important fact we will use is the Rademacher theorem,
namely a Lipschitz function $f$ on $X$ is differentiable $\vol$-a.e..
This easily follows from the usual Rademacher theorem for $f \circ \phi^{-1}$
(see \cite[Corollary~2.14]{Bert}). It follows from $(\ref{eq:Taylor})$ that
$\sqrt{\langle \nabla f(x),\nabla f(x) \rangle}$ coincides with
the local Lipschitz constant $\gl{f}(x)$ given by 
\[ \gl{f}(x) :=\lims_{y \to x} \frac{|f(x)-f(y)|}{d(x,y)}. \]

Based on the notion of gradient vector, 
we define the Sobolev space and the Dirichlet energy as follows
(see \cite{Kuw_Mac_Shi} and \cite{Kuw_Shi03} for details).
For a function $f:X \to \R$ such that
$f \circ \phi^{-1} \in W^{1,2}(\phi(U_{\phi}))$ for all $\phi \in \Phi$,
we introduce the weak gradient vector $\nabla f(x) \in K_x$
for a.e.\ $x \in V_{\phi}$ as the element corresponding to the weak
gradient vector $\nabla(f \circ \phi^{-1})(\phi(x))$.
We define the Sobolev space $W^{1,2} (X)$ and
the Dirichlet energy $\mathcal{E}$ by 
\begin{align*}
W^{1,2} (X) & : = 
\left\{ 
f \in L^2 ( X, \vol ) 
  \,\left| \, 
\int_X \langle \nabla f,\nabla f \rangle \,d\vol<\infty 
\right. \right\}, 
\\
\mathcal{E}(f,g) & := 
\int_X \langle \nabla f,\nabla g \rangle \,d\vol \qquad
{\rm for}\ f,g \in W^{1,2}(X).
\end{align*}
We do not divide $\mathcal{E}$ by $2$ for notational simplicity.
Note that $\mathcal{E}$ coincides with the energy functional introduced by
Korevaar and Schoen~\cite{Kor-Sch}
(we can reduce the argument to the Euclidean case by using a $(1+\ep)$-bi-Lipschitz chart;
see \cite[Theorem~6.2]{Kuw_Shi03}).
We also remark that the set of Lipschitz functions
$C^{\mathrm{Lip}}(X)$ is dense in $W^{1,2}(X)$
with respect to the Sobolev norm
$\|f\|_{W^{1,2}}^2=\|f\|_{L^2}^2+\mathcal{E}(f, f)$
(\cite[Theorem~1.1]{Kuw_Mac_Shi}). 

Furthermore, if $f$ is a Lipschitz function, then the weak gradient vector
$\nabla f(x)$ coincides with the gradient vector as in $(\ref{eq:Taylor})$ a.e.\ $x$,
and hence $\sqrt{\langle \nabla f(x),\nabla f(x) \rangle}=\gl{f}(x)$
holds a.e.\ $x$.
Therefore $\mathcal{E}$ also coincides with Cheeger's energy functional
(\cite{Chee99}), because the local Lipschitz constant is the minimal generalized
upper gradient (\cite[Theorem~6.1]{Chee99}) and Lipschitz functions are
dense in both Sobolev spaces (thanks to the weak Poincar\'e inequality 
for upper gradients and the volume doubling condition, 
\cite[Theorem~4.24]{Chee99}). 
Indeed, in our framework, the volume doubling condition directly follows from 
the Bishop-Gromov volume comparison theorem and
the Poincar\'e inequality is a consequence 
of \cite{Kuw_Shi01} and \cite{Ranj_PI}, for instance.
\bigskip

By following the general theory of bilinear forms, 
there exists a nonpositive
selfadjoint operator 
$( \lap , D ( \lap ) )$ on $L^2 ( X , \vol )$ 
associated with $(\mathcal{E} , W^{1,2} (X) )$. 
It is characterized by the following identity: 
\begin{equation}
\label{eq:bypart}
\mathcal{E} ( g , f ) 
= 
- \int_X g \lap f  \, d \vol  
, \quad  
f \in D ( \lap ), \, g \in W^{1,2} (X)
. 
\end{equation}
We call $\lap$ the Laplacian as in the classical case. 
Based on a general theory of functional analysis, 
the one-parameter semigroup of 
contractive symmetric linear operators 
$T_t = \e^{t \lap}$ on $L^2 ( X, \vol )$ 
is defined associated with $\lap$.
For any $f \in L^2 ( X , \vol )$, 
$T_t f$ solves the (linear) heat equation 
$\partial_t u = \lap u$ 
with $u (0 , \cdot ) = f$ 
in the sense that 
$T_t f \in D (\lap)$ for $t > 0$ 
and 
\begin{align*}
\lim_{t \downarrow 0} 
\frac{T_{t} f - f}{t} 
& = 
\lap f 
\quad 
\mbox{
  in $L^2 ( X, \vol )$ 
  for $f \in D ( \lap )$, 
} 
\\
\lim_{t \downarrow 0} 
T_t f 
& = 
f 
\quad 
\mbox{ 
  in $L^2 ( X, \vol )$ 
  for $f \in L^2 ( X, \vol)$. 
} 
\end{align*} 
Note that $T_t$ is Markovian in the sense 
that $0 \le T_t f \le 1$ holds 
whenever $0 \le f \le 1$. 
As shown in \cite[Theorem~1.5]{Kuw_Mac_Shi}, 
there exists a continuous function 
$(t, x, y ) \mapsto p_t (x,y)$ 
on $(0,\infty) \times X \times X$ 
satisfying the following properties: 
\begin{enumerate}
\item
For any $f \in L^2 ( X, \vol )$, 
$t > 0$ and $\vol$-a.e.\ $x \in X$,   
\begin{equation} \label{eq:D}
T_t f (x) 
= 
\int_X p_t ( x, y ) f (y) \,\vol (dy) .
\end{equation}
\item
For any $s,t > 0$ and $x,y \in X$, 
\begin{align}
\label{eq:S}
p_t (x,y) 
& =
p_t (y,x) , 
\\
\label{eq:CK}
p_{s+t} (x,y) 
& =
\int_X p_s (x,z) p_t (z,y) \,\vol (dz) ,
\\
\label{eq:C}
\int_X p_t (x,z) \,\vol (dz) 
& = 1 ,
\\ \nonumber
p_t (x,y) 
& > 0 .
\end{align} 
\end{enumerate}
See \cite[Theorems~1.4, 1.5(3)]{Kuw_Mac_Shi} 
for the continuity of $p_t (x,y)$. 
The equality \eqref{eq:C} follows from 
the fact that $1 \in W^{1,2} (X)$ and $\mathcal{E} (1,1) = 0$,
because $X$ is assumed to be compact. 
As in the classical case, we call 
$p_t (x,y)$ the heat kernel. 
The existence and these properties of $p_t (x,y)$ are 
deduced from the Poincar\'e inequality for $\mathcal{E}$ 
and the volume doubling condition, 
together with results in \cite{Sturm_Harnack}. 
\begin{remark}
The bilinear form $( \mathcal{E} , W^{1,2} (X) )$ is 
a symmetric strongly local regular Dirichlet form 
(see \cite{Kuw_Mac_Shi} for it and further details; 
see \cite{BH, FOT} for basics on Dirichlet forms). 
Moreover, $T_t$ enjoys the strong Feller property. 
As a result, 
there exists a diffusion process 
$( ( X_t )_{ t \ge 0} , ( \mathbb{P}_x )_{x \in X} )$ 
on the whole space $X$ 
associated with $( \mathcal{E}, W^{1,2} (X) )$ 
in the sense that $\mathbb{E}_x [ f ( X_t ) ] = T_t f (x)$ 
for $f \in C(X)$, $x \in X$ and $t > 0$. 
\end{remark} 

Denote by $\prob{X}$ the set of all Borel probability measures on $X$.
Let us define a positive Borel measure 
$T_t \nu$ for $\nu \in \prob{X}$ and $t \ge 0$ by 
\[
T_t \nu ( dy ) 
:= 
\begin{cases}
\displaystyle
\abra{ 
  \int_X p_t ( x , y ) \,\nu (dx) 
} 
\vol (dy)
 
& t > 0 ,
\\
\nu ( dy ) 
& t = 0 .
\end{cases}
\]
Thanks to \eqref{eq:C}, $T_t \nu \in \prob{X}$ holds. 
By definition, 
$T_t \nu$ is absolutely continuous 
with respect to $\vol$ for $t > 0$. 
When $d\nu = f d\vol$, 
it holds $d T_t \nu = T_t f d \vol$. 
In this paper, we call the evolution $(t,\nu) \mapsto T_t \nu$
\emph{the gradient flow of the Dirichlet energy} 
(since the Dirichlet energy is a functional 
on the $L^2$-space of functions, 
this terminology should be interpreted in an extended sense). 
Indeed, it is easy to see from
\[ \frac{1}{2}\mathcal{E}(f+\ep g,f+\ep g)
 =\frac{1}{2}\mathcal{E}(f,f) -\ep \int_X g\lap f \,d\vol+O(\ep^2) \]
that the Radon-Nikodym derivative $d T_t \nu / d \vol$ 
is the gradient flow of $\mathcal{E}/2$
with respect to the $L^2$-norm.

Before closing this subsection, 
we review the derivation property of 
$W^{1,2} (X)$. 
It is formulated as follows: 
For $f_1 , \ldots , f_k , g \in W^{1,2} (X) \cap L^\infty ( X, \vol )$ 
and $\Phi \: : \: \R^ k \to \R$ which is $C^1$ 
on the range of $(f_1 , \ldots , f_k )$, 
$\Phi ( f_1 , \ldots , f_k )$ belongs to $W^{1,2} (X)$ 
and 
\begin{equation}
\label{eq:chain}
\gss{\Phi ( f_1 , \ldots , f_k )}{g} 
= 
\sum_{j=1}^{k}
\frac{\partial \Phi}{\partial x_j} ( f_1 , \ldots , f_{k} ) 
\gss{f_j}{g} 
\quad 
\mbox{$\vol$-a.e..}
\end{equation}
This identity directly follows from the definition of $W^{1,2} (X)$ or 
from the strong locality of the Dirichlet form 
$(\mathcal{E} , W^{1,2} (X) )$ 
(see \cite[Corollary~I.6.1.3]{BH} and \cite[Section 3.2]{FOT} for the latter).
\end{subsection}

\begin{subsection}{Gradient flows in the Wasserstein space}

We next introduce the Wasserstein space
and a purely metric notion of gradient flows in it.
We refer to \cite{AGS} and \cite{Villani08} for the basic theory
as well as the recent diverse developments.

Given $\mu,\nu \in \prob{X}$, a probability measure
$\pi \in \prob{X \times X}$ is called a coupling of $\mu$ and $\nu$
if $\pi(A \times X)=\mu(A)$ and $\pi(X \times A)=\nu(A)$ hold
for all Borel sets $A \subset X$.
Then, for $1 \le p < \infty$, 
we define the $L^p$-Wasserstein distance as
\[ W_p(\mu,\nu):=\inf_{\pi} \bigg( \int_{X \times X}
 d(x,y)^p \,\pi(dxdy) \bigg)^{1/p}, \]
where $\pi$ runs over all couplings of $\mu$ and $\nu$.
In most parts, we work in the quadratic case $p=2$. 
The $L^2$-Wasserstein space $(\prob{X},W_2)$ 
becomes a metric space and inherits 
the compactness from $(X,d)$.
Moreover, $(\prob{X},W_2)$ is a geodesic space.
If $\mu$ is absolutely continuous with respect to $\vol$,
then a minimal geodesic $(\mu_t)_{t \in [0,1]}$ from $\mu$
to any $\nu$ is unique and $\mu_t$ is also absolutely continuous
for all $t \in (0,1)$ (see \cite{Bert} for a more detailed
characterization of $\mu_t$, and \cite{FJ} for the absolute continuity).

For $\mu \in \prob{X}$, we define the relative entropy by
\[ \Ent(\mu):=\int_X \rho \log \rho \,d\vol \]
when $\mu=\rho \vol$ with $\rho \in L^1 ( X, \vol)$, 
and $\Ent(\mu):=\infty$ otherwise.
Set $\pacs{X}:=\{ \mu \in \prob{X} \,|\, \Ent(\mu)<\infty \}$.
Note that $\Ent$ is lower semi-continuous with respect to $W_2$
and satisfies $\Ent(\mu) \ge -\log\vol(X)$ by Jensen's inequality.

\begin{definition}[The curvature-dimension condition]\label{df:CD}
For $K \in \R$, we say that $(X,d,\vol)$ satisfies
the \emph{curvature-dimension condition} $\CD(K,\infty)$ if $\Ent$ is
$K$-geodesically convex in the sense that any pair $\mu,\nu \in \prob{X}$ admits
a minimal geodesic $(\mu_t)_{t \in [0,1]}$ from $\mu$ to $\nu$ such that
\[ \Ent(\mu_t) \le (1-t)\Ent(\mu) +t\Ent(\nu) -\frac{K}{2}(1-t)tW_2(\mu,\nu)^2 \]
holds for all $t \in [0,1]$.
\end{definition}

We remark that the above inequality is obvious
if $\mu \not\in \pacs{X}$ or $\nu \not\in \pacs{X}$.
Therefore it is sufficient to consider $\mu,\nu \in \pacs{X}$,
and then a minimal geodesic between them is unique.

\begin{remark} \label{rem:CD}
(i) The curvature-dimension condition $\CD(K,\infty)$ for general metric measure
spaces is introduced and studied independently in \cite{SturmI} and \cite{LV-CD},
and known to be equivalent to the lower Ricci curvature bound $\Ric \ge K$
for complete Riemannian manifolds equipped with the Riemannian distance and the volume measure (\cite{vRS}).

(ii) The condition $\CD(K,N)$ for $N \in (1,\infty)$ is also introduced in
\cite{SturmII} and \cite{LV-CD2}.
In general, $\CD(K,N)$ implies $\CD(K,\infty)$.
In the Riemannian case, $\CD(K,N)$ is equivalent to $\Ric \ge K$ and $\dim \le N$.

(iii) It is recently demonstrated in \cite{PetCD} and \cite{ZZ1} that
$n$-dimensional Alexandrov spaces of curvature bounded from below by $k$ satisfy
$\CD((n-1)k,n)$ (and hence $\CD((n-1)k,\infty)$), as is naturally expected
from the relation between the sectional and the Ricci curvatures.
\end{remark}

There is a well established theory on the gradient flow of geodesically
convex functions as comprehensively discussed in \cite{AGS}.
For later convenience, we recall a couple of notions in a general form.
We say that a curve $(\mu_t)_{t \in I} \subset \prob{X}$ on an interval
$I \subset \R$ is absolutely continuous if there is $f \in L^1(I)$ such that
\begin{equation}\label{eq:accurve}
W_2(\mu_t,\mu_s) \le \int_t^s f(r) \,dr
\end{equation}
for all $t,s\in I$ with $t \le s$.
Note that absolutely continuous curves are continuous.
For an absolutely continuous curve $(\mu_t)_{t \in I}$,
the metric derivative
\[ |\dot{\mu}_t|:=\lim_{h \to 0}\frac{W_2(\mu_t,\mu_{t+h})}{|h|} \]
is well-defined for a.e.\ $t \in I$ (\cite[Theorem~1.1.2]{AGS}).
Moreover, $|\dot{\mu}_t|$ belongs to $L^1(I)$ and is the minimal function
for which \eqref{eq:accurve} holds.
Given a functional $E:\prob{X} \to \R \cup \{+\infty\}$,
we consider a gradient flow of $E$ solving
``$\dot{\mu}_t=-\nabla E(\mu_t)$'' in the following sense.
For $\mu \in \prob{X}$ with $E (\mu) < \infty$,
we define the local slope as
\begin{equation}\label{eq:slope}
\ls{E}(\mu) :=\lims_{\nu \to \mu}
 \frac{\max\{E(\mu)-E(\nu),0\}}{W_2(\mu,\nu)}.
\end{equation}
If $E$ is $K$-geodesically convex, then we have
\begin{equation}\label{eq:strongug}
|E(\mu_t)-E(\mu_s)| \le \int_t^s
 |\dot{\mu}_r| \cdot \ls{E}(\mu_r) \,dr
\end{equation}
for all $t,s \in I$ with $t<s$ along any absolutely continuous curve
$(\mu_t)_{t \in I}$ with values in $\prob{X}$ with $E (\mu_t ) < \infty$. 
As a consequence, it holds
\begin{equation}\label{eq:fromug}
E(\mu_t) \le E(\mu_s)+\frac{1}{2}\int_t^s|\dot{\mu}_r|^2 \,dr
 +\frac{1}{2}\int_t^s \ls{E}^2(\mu_r) \,dr
\end{equation}
for all $t<s$.
From $(\ref{eq:strongug})$ and $(\ref{eq:fromug})$,
it is natural to give the following definition.

\begin{definition}[Gradient flows]\label{df:GF}
Let $E:\prob{X} \to \R \cup \{+\infty\}$ be a $K$-geodesically convex functional.
We say that an absolutely continuous curve
$(\mu_t)_{t \in [0,\infty)}$ in $\prob{X}$ 
is a \emph{gradient flow} of $E$ provided 
$E ( \mu_t ) < \infty$ for $t \ge 0$ and 
\begin{equation}\label{eq:defgf}
E(\mu_t)= E(\mu_s)+\frac12 \int_t^s|\dot{\mu}_r|^2 \,dr
 +\frac12 \int_t^s \ls{E}^2(\mu_r) \,dr
\end{equation}
for all $0 \le t<s$.
The equation \eqref{eq:defgf} is called the \emph{energy dissipation identity}.
\end{definition}

The existence of such a gradient flow comes from the general theory presented
in \cite[Corollary~2.4.11]{AGS}.
Furthermore, as shown in \cite[Theorem~6.2]{Ohta_grad-Alex}
(see also \cite[Theorem~4.2]{Gigli-Ohta10}),
the gradient flow produces a contraction semigroup in the sense that
for any $\mu,\nu \in \prob{X}$ with $E(\mu) < \infty$ and $E (\nu) < \infty$,
the gradient flows $(\mu_t)_{t \in [0,\infty)}, (\nu_t)_{t \in [0,\infty)}$
starting from $\mu,\nu$ satisfy
\begin{equation}\label{eq:contr}
W_2(\mu_t,\nu_t)\leq \e^{-Kt}W_2(\mu,\nu)
 \qquad {\rm for\ all}\ t\geq 0,
\end{equation}
where $K$ is the modulus of convexity of $E$.
In particular, the uniqueness follows from considering $\mu=\nu$.
(Though the strategy of the construction in \cite{Ohta_grad-Alex} is different from
\cite{AGS}, the resulting flow is the same by uniqueness,
see \cite[Remark~2.7]{Gigli-Ohta10}.)

\begin{remark}\label{rem:GF}
(i) In \cite[Theorem~7]{Sav} the contractivity (from the geodesical convexity)
is shown on spaces satisfying the \emph{local angle condition}.
Alexandrov spaces satisfy this condition.

(ii) The first author \cite[Theorem~15]{heatmm} proved the uniqueness
of the gradient flow of the relative entropy on general metric measure spaces
satisfying $\CD(K,\infty)$, without relying on the contractivity.
In this generality, the contractivity fails. 
More precisely, 
the heat flow on a finite-dimensional (Minkowski) normed space
is not contractive (except for inner product spaces; see \cite{OS2}).
\end{remark}

In the sequel, we mainly study the gradient flow of the relative entropy.
Then $(\ref{eq:contr})$ allows us to continuously and uniquely
extend the gradient flow semigroup to the full $\prob X$
(since the closure of $\pacs{X}$ is $\prob{X}$).
Such an extension also satisfies $(\ref{eq:contr})$.
\end{subsection}
\end{section}

\begin{section}{Identification of the two gradient flows}\label{sc:main}

This section contains our main result:

\begin{theorem}[Identification of the two gradient flows]
\label{thm:main}
Let $(X,d)$ be a compact $n$-dimensional Alexandrov space without boundary.
For any $\nu \in \prob{X}$, the gradient flow of the Dirichlet energy starting from $\nu$
is the gradient flow of the relative entropy, and vice versa.
\end{theorem}

Recall that $\CD(K,\infty)$ holds with $K=(n-1)k$, and we use this,
e.g., in the proof of Proposition~\ref{prop:slope} below.
The main technical obstacle in the proof of this theorem is to let the $L^2$
and the $W_2$ structures `interact'.
Our strategy consists in picking a gradient flow of the Dirichlet energy
and in proving that it obeys the energy dissipation identity
$(\ref{eq:defgf})$ for the relative entropy.
We start with a bound on the local slope $(\ref{eq:slope})$
by the Fisher information.

\begin{proposition}[Bound on the local slope]\label{prop:slope}
Let $\mu = f \vol \in \pacs{X}$ with $f \in W^{1,2} (X)$.
Then it holds
\[ \ls{\Ent}^2( \mu ) \leq \int_X \frac{\gss f f}{f} \,d\vol. \]
\end{proposition}

\begin{Proof}
We first assume that $f$ is Lipschitz and bounded away from 0.
Then we know that $\gss f f=\gl{f}^2$ $\vol$-a.e.,
so that the conclusion follows from \cite[Theorem~20.1]{Villani08} 
together with $\CD(K,\infty)$.
Thus all we need to do is to proceed by approximation.
Suppose that $0<c\leq f$ $\vol$-a.e.\ for some $c \in \R$.
Since $C^{\mathrm{Lip}} (X)$ is dense in $W^{1,2} (X)$ and $\sqrt{f} \in W^{1,2} (X)$, 
we can find a sequence $\{g_i\}_{i \in \N}$
of Lipschitz functions such that $g_i$ converges to $\sqrt{f}$
as $i \to \infty$ with respect to the Sobolev norm.
Substituting $C_i \max\{ g_i,\sqrt{c} \}$ for some $C_i>0$ if necessary,
we can assume that $\sqrt{c_i} \leq g_i$ for some $c_i>0$
$\vol$-a.e.\ as well as
$\| g_i \|_{L^2} = 1$ for all $i$. 
Set $f_i:=g_i^2$. 
As $f_i$ is Lipschitz and bounded away from $0$, we have
\[ \ls{\Ent}^2(f_i\vol) \leq \int_X \frac{\gss {f_i}{f_i}}{f_i} \,d\vol. \]
On the one hand, the right-hand side is equal to $4\mathcal E(g_i,g_i)$
and converges to
$4\mathcal E(\sqrt f,\sqrt f)=\int_X {\gss f f}/f \,d\vol$ by construction.
On the other hand, we deduce from \cite[Corollary~2.4.10]{AGS} that
\begin{equation} \label{eq:lsc}
 \ls{\Ent}(f\vol) \leq \limi_{i\to\infty}\ls{\Ent}(f_i\vol). 
\end{equation}
Combining these, we prove the claim for $f$.
It remains to remove the assumption that $f$ is bounded away from $0$.
To do this we just consider $f_i :=(1+i^{-1})^{-1} ( f + i^{-1} )$ 
and apply \eqref{eq:lsc} again.
\end{Proof}

Let us turn to considering the gradient flow of the Dirichlet energy $(T_t \nu)_{t \in [0,\infty)}$. 
For simplicity of notations, we write $f_t = d T_t \nu / d \vol$
in the remainder of the section.
By \eqref{eq:D}, $f_t \in C(X) \subset L^2 ( X , \vol )$ holds for $t > 0$. 
Then \eqref{eq:CK}, \eqref{eq:S} and \eqref{eq:D} imply 
$f_t \in D ( \lap )$. For any $t >0$,  there is  $\ep_t > 0$ such that $f_t \ge \ep_t$ holds 
since $p_t (x,y)$ is positive and continuous. 
As a well-known fact, the following bound also holds. 

\begin{lemma}[Maximum principle]
\label{le:max}
Let $f_0 \in L^2 ( X, \vol )$. 
If $f_0\geq c$ a.e., then $T_t f_0 \geq c$ a.e. for every $t\geq 0$.
The same holds for bounds from above.
\end{lemma}

\begin{Proof} 
Take $A \subset X$ measurable. 
Since $T_t$ is Markovian, $T_t 1_A \ge 0$ $\vol$-a.e.. 
Thus 
\begin{equation*} 
\int_A ( T_t f_0 - c ) \,d \vol
  =
\int_A T_t ( f_0 - c ) \,d \vol 
  = 
\int_X T_t 1_A \cdot ( f_0 - c ) \,d \vol 
\ge 0 . 
\end{equation*} 
Since $A$ is arbitrary, the assertion holds. 
Bounds from above follow by applying the same argument to $-f_0$. 
\end{Proof}
By virtue of Lemma~\ref{le:max} with the remark before it, 
for any $\ep>0$, there are $c, C >0$ such that 
\begin{equation} \label{eq:PIP}
c\leq f_t \leq C \quad \mbox{for all $t\geq \ep$.}
\end{equation}
\begin{proposition}[Entropy dissipation]\label{prop:entdiss}
Let $\nu \in \prob{X}$ and $f_t=dT_t \nu/d\vol$.
Then the function $t\mapsto\Ent(f_t\vol)$ is locally Lipschitz in $(0,\infty)$ and,
moreover, it holds
\[ \frac d{dt}\Ent(f_t\vol)=-\int_X\frac{\gss{f_t}{f_t}}{f_t} \,d\vol
 \qquad {\it a.e.}\ t. \]
\end{proposition}

\begin{Proof}
As the function $e(s):=s\log s$ is $C^1$ in $[c, C]$ and
$t\mapsto f_t$ is locally Lipschitz in $(0,\infty)$ with values in $L^2(X, \vol)$,
we deduce from \eqref{eq:PIP} that 
$t\mapsto \Ent(f_t\vol)$ is locally Lipschitz in $(0,\infty)$.
Applying formulas \eqref{eq:bypart} and \eqref{eq:chain}
and recalling $\int_X \lap f_t \,d\vol=0$, we obtain
\begin{align*}
&\frac d{dt}\Ent(f_t\vol)
 =\int_X e'(f_t)\lap f_t \,d\vol =\int_X \big(\log(f_t)+1\big)\lap f_t \,d\vol \\
&=-\int_X \gss{\log(f_t)}{f_t} \,d\vol =-\int_X \frac{\gss{f_t}{f_t}}{f_t} \,d\vol.
\end{align*}
\end{Proof}

For the next argument, 
we briefly recall some properties of the Hamilton-Jacobi semigroup 
in our context. 
For $f \in C^{\mathrm{Lip}} (X)$ and $t>0$, 
we define $Q_tf:X\to\R$ by 
\begin{equation}
\label{eq:hopflax}
Q_tf(x):=\inf_{y\in X} \left[ f(y)+\frac{d^2(x,y)}{2t} \right].
\end{equation}
Also, set $Q_0f:=f$. 
The following is shown in 
\cite[Theorem~2.5(iv)]{BEHM_HJ} and \cite[Theorem~2.5(viii)]{HJ}
in the framework of general metric measure spaces supporting 
the volume doubling condition 
and the Poincar\'e inequality for upper gradients. 

\begin{proposition}[Hamilton-Jacobi semigroup]
\label{prop:HJ}
It holds $Q_tf \in C^{\mathrm{Lip}} (X)$ for every $t\geq 0$, 
the map $[0,\infty)\ni t\mapsto Q_tf\in C(X)$ is Lipschitz
in the uniform norm, and 
\begin{equation}
\label{eq:HJ}
\frac {d}{dt} Q_tf(x)+\frac{\gl{Q_tf}^2(x)}{2}=0,
\end{equation}
for almost every $t$, $x$.
\end{proposition}

\begin{remark}
The equation \eqref{eq:hopflax} has been called the Hopf-Lax formula
or the Moreau-Yosida approximation also in the literature. 
The former name is mainly used in the PDE context 
to a \emph{special} solution to the Hamilton-Jacobi equation \eqref{eq:HJ}. 
The latter one is typically used 
for an approximation of an (unbounded) operator in functional analysis. 
\end{remark}

\begin{proposition}[Absolute continuity with respect to $W_2$]
\label{prop:acw2}
For $\nu \in \prob{X}$, the curve $t\mapsto T_t \nu$ is absolutely continuous 
in the space $(\prob X,W_2)$ and 
its metric speed $|\dot{T_t\nu}|$ satisfies 
\begin{equation}
\label{eq:boundspeed}
|\dot{T_t \nu}|^2 \leq \int_X \frac{\gss{f_t}{f_t}}{f_t} \,d\vol,
\qquad {\it for\ a.e.}\ t.
\end{equation}
\end{proposition}

\begin{Proof}
Fix $t,s>0$. By the Kantorovich duality (cf., e.g., \cite[Theorem~6.1.1]{AGS} and \cite[Theorem~5.10]{Villani08})
together with \eqref{eq:hopflax}, we obtain 
\begin{equation} \label{eq:step1} 
\frac12 W_2^2\big( T_{t}\nu , T_{t+s}\nu \big)
=
\sup_{\varphi \in C^{\mathrm{Lip}} (X) } 
\left[ 
    \int_X (Q_1 \varphi) f_{t+s} \,d \vol - \int_X \varphi f_{t} \,d \vol
\right] . 
\end{equation} 
By Proposition~\ref{prop:HJ}, 
the map $r \mapsto Q_{r} \varphi$ from $[0,1]$ to $L^2 ( X, \vol )$ is Lipschitz. 
Moreover, it is differentiable in $L^2 (X , \vol)$ 
and the derivative is determined by \eqref{eq:HJ}. 
The curve $[0,1]\ni r\mapsto f_{t+rs}\in L^2(X,\vol)$ is Lipschitz as well.
Thus $[0,1]\ni r\mapsto (Q_{r}\varphi) f_{t+rs}\in L^1 (X,\vol)$ is Lipschitz and its derivative can be calculated with the Leibniz rule.
Thus we have 
\begin{align} \nonumber
\int_X (Q_1 \varphi) & f_{t+s} \,d\vol-\int_X \varphi f_t \,d\vol\\
\nonumber
&=\int_0^1\frac d{dr} \bigg[ \int_X (Q_{r} \varphi) f_{t+rs} \,d\vol \bigg] \,dr\\
\label{eq:step2}
&=\int_0^1\int_X \left(-\frac{\gl{Q_{r} \varphi}^2}{2}f_{t+rs}
 +s(Q_{r} \varphi) \lap f_{t+rs} \right) d\vol\,dr . 
\end{align}
Using formulas \eqref{eq:bypart}, \eqref{eq:PIP} and the trivial inequality
\[
-\gss g{\tilde g}\leq \frac1{2s}\gss g g+\frac{s}{2}\gss{\tilde g}{\tilde g},\qquad \mbox{$\vol$-a.e. for $s > 0$},
\] 
we have
\[
\begin{split}
\int_X (Q_{r}& \varphi) \lap f_{t+rs} \,d\vol\\
&=-\int_X \gss{Q_{r} \varphi }{f_{t+rs}} d\vol\\
&=-\int_X \la \nabla Q_{r} \varphi , \frac{\nabla f_{t+rs}}{f_{t+rs}} \ra f_{t+rs} \,d\vol\\
&\leq \frac{1}{2s} \int_X \gss{Q_{r} \varphi}{Q_{r} \varphi} f_{t+rs} \,d\vol  
+ \frac{s}{2} \int_X \frac{\gss{f_{t+rs}}{f_{t+rs}}}{f_{t+rs}} \,d\vol.
\end{split}
\]
Plugging this inequality in \eqref{eq:step2} and recalling that
$\gss{Q_{r} \varphi}{Q_{r}\varphi}=\gl{Q_{r}\varphi}^2$ $\vol$-a.e.\
(since $Q_{r} \varphi$ is Lipschitz), we obtain
\[
\begin{split}
\int_X (Q_1\varphi) f_{t+s} \,d\vol-\int_X \varphi f_t \,d\vol 
\leq 
\frac{s^2}{2} \int_0^1\int_X \frac{\gss{f_{t+rs}}{f_{t+rs}}}{f_{t+rs}} \,d\vol\,dr.
\end{split}
\]
This bound does not depend on $\varphi$, thus from \eqref{eq:step1} we deduce
\begin{equation}
\label{eq:step3}
W_2^2\big(T_{t}\nu, T_{t+s}\nu \big)\leq
 {s^2}\int_0^1\int_X \frac{\gss{f_{t+rs}}{f_{t+rs}}}{f_{t+rs}} \,d\vol\,dr.
\end{equation}
Since we have \eqref{eq:PIP} and 
the Dirichlet energy decreases along the flow $t\mapsto f_t$, we obtain 
\[
W_2^2\big( T_{t}\nu ,T_{t+s}\nu \big)\leq
 \frac {s^2}{c}\int_0^1\int_X {\gss{f_{t+rs}}{f_{t+rs}}} \,d\vol\,dr\leq
 \frac {s^2}{c}\int_X {\gss{f_{t}}{f_{t}}} \,d\vol,
\]
which gives that the map $t\mapsto T_{t} \nu \in \prob X$ is locally Lipschitz. The bound \eqref{eq:boundspeed} follows directly from \eqref{eq:step3}
(by recalling the definition of the absolutely continuous curves $(\ref{eq:accurve})$).
\end{Proof}

Now we are ready to prove Theorem~\ref{thm:main}.

\begin{Proofof}\emph{Theorem}~\ref{thm:main}.
By Propositions~\ref{prop:acw2}, \ref{prop:slope} and \ref{prop:entdiss},
we have
\[ \frac12|\dot{f_t\vol}|^2+\frac12\ls{\Ent}^2(f_t\vol)
 \leq \int_X\frac{\gss{f_t}{f_t}}{f_t} \,d\vol
 =-\frac d{dt}\Ent(f_t\vol) \]
a.e.\ $t$.
As the reverse inequality $(\ref{eq:fromug})$ is always true,
equality $(\ref{eq:defgf})$ holds for all $0 < t \le s$. 
Since $(f_t\vol)_{t \in [0,\infty)}$ is continuous also at $t = 0$, 
it is the gradient flow of the relative entropy.
The converse immediately follows from the uniqueness of the gradient flow
of the relative entropy.
\end{Proofof}

\begin{remark} \label{rem:initial}
(i) When $\nu = f_0 \vol$ with $f_0 \in L^2 (X, \vol)$, 
we can give a proof of Theorem~\ref{thm:main} 
without relying on the positivity improving property \eqref{eq:PIP}. 
Indeed, by virtue of the contraction property of both flows, 
it is possible to use Propositions~\ref{prop:entdiss}, \ref{prop:acw2}
only for $0<c\le f_0 \le C$ and prove the theorem via approximation.
It suggests that our argument possibly works in a more general framework 
where the existence of the density $p_t (x,y)$ does not follow from
the theory of Dirichlet forms.
In such a case, we could `construct' $p_t(x,\cdot)$ as the gradient flow
starting from the Dirac measure $\delta_x \in \prob{X}$
(via the contraction property).

(ii) In the proof of Theorem~\ref{thm:main}, 
we heavily rely on the fact 
$\gss{f}{f} = \gl{f}^2$ for $f \in C^{\mathrm{Lip}} (X)$ 
for which we used the local structure of Alexandrov spaces. 
Actually, we use $\gl{f}$ in Propositions~\ref{prop:slope}, \ref{prop:acw2} 
and $\gs{f}$ in Propositions~\ref{prop:entdiss}, \ref{prop:acw2}.  
\end{remark}
\end{section}

\begin{section}{Applications}\label{sc:appl}

In this section, we assume that $(X,d)$ is a compact Alexandrov space
without boundary satisfying $\CD(K,\infty)$, and 
prove some applications of Theorem~\ref{thm:main}. 
It should be stressed that the `sectional curvature bound' $k$ 
in the sense of Alexandrov appears nowhere in the sequel,
and the `Ricci curvature bound' $K$ is essential instead
(recall Remark~\ref{rem:CD}).
Indeed, those results involving $K$ are 
natural extensions of the corresponding ones 
on a Riemannian manifold with $\Ric \ge K$. 
In this sense, the emergence of $K$ instead of $k$ 
is natural and gives sharper estimates. 

Since the gradient flow of the Dirichlet energy $T_t \nu$ is obviously 
linear and symmetric, we immediately obtain the following: 
\begin{theorem}[Linearity and symmetry] 
For $\nu \in \prob{X}$, 
let $( \mu_t^\nu )_{t \ge 0}$ be 
the gradient flow of the relative entropy 
on $\mathscr{P} (X)$ with $\mu_0^\nu = \nu$. 
Then the following hold true. 
\begin{enumerate}
\item
For 
$\nu_0 , \nu_1 \in \prob{X}$, 
$\lambda \in [0,1]$ and 
$t \ge 0$, 
\begin{equation*} 
\mu_t^{ ( 1 - \lambda ) \nu_0 + \lambda \nu_1 }
  =  
( 1 - \lambda ) \mu_t^{\nu_0}
  + 
\lambda \mu_t^{\nu_1} .
\end{equation*}
\item
For $f,g \in L^1 (X , \vol)$ 
with $f ,g \ge 0$ and 
$\| f \|_{L^1} = \| g \|_{L^1} = 1$, 
\begin{equation*} 
\int_X f \, d \mu_t^{g \vol} 
  = 
\int_X g \, d \mu_t^{f \vol} . 
\end{equation*}
\end{enumerate}
\end{theorem} 

We remark that the linearity, but not symmetry, also follows from the gradient flow approach
under the local angle condition (\cite[Theorem~8]{Sav}).
In general, these properties are completely nontrivial,
and the linearity indeed fails in the Finsler setting (\cite{heatfinsl}).

A new property for $( T_t \nu )_{t \ge 0 , \nu \in \prob{X} }$ 
coming from our identification with the gradient flow of the relative entropy is 
the $L^2$-Wasserstein contraction \eqref{eq:contr}.
Together with \cite[Corollary~3.4]{K9}, we obtain the following: 
\begin{theorem}[Contraction for the heat flow] 
\label{thm:contraction}
For $1 \le p \le 2$, 
\begin{equation} \label{eq:contraction2}
W_p ( T_t \nu_0 , T_t \nu_1 ) 
\le 
\e^{- K t} W_p ( \nu_0 , \nu_1 ) 
\end{equation}
holds 
for every $\nu_0 , \nu_1 \in \prob{X}$. 
\end{theorem}
Furthermore, by the duality result \cite[Theorem~2.2]{K9}, 
Theorem~\ref{thm:contraction} yields 
the following Bakry-\'{E}mery type $L^2$-gradient estimate 
\begin{equation} \label{eq:BE}
\gl{T_t f} (x) \le \e^{-Kt} T_t ( \gl{f}^2 ) (x)^{1/2} 
\end{equation}
for any $f \in C^{\mathrm{Lip}} (X)$ and $x \in X$. 
By combining \eqref{eq:BE} 
with the regularity of the heat kernel, 
we can extend it to $f \in W^{1,2} (X)$ as follows: 
\begin{theorem}[Gradient estimate] 
\label{thm:BE}
Let $f \in W^{1,2} (X)$ and $t > 0$. 
Then $T_t f \in C^{\mathrm{Lip}} (X)$ 
and 
\begin{equation} \label{eq:BE2}
\gl{T_t f} (x) \le \e^{-Kt} T_t ( \gs{f}^2 ) (x)^{1/2} 
\end{equation}
holds for all $x \in X$. 
In particular, 
\begin{equation} \label{eq:BE3}
\gs{T_t f} (x) \le \e^{-Kt} T_t ( \gs{f}^2 ) (x)^{1/2} 
\quad 
{\it for\ a.e.}\ x
\end{equation}
and 
$\gl{T_t f} \le 
\e^{-Kt} \sqrt{\| T_t \|_{L^1 \to L^\infty} \mathcal{E} (f,f)}$ 
hold. 
\end{theorem}
\begin{Proof}
Take $\{ f_i \}_{i \in \N} \subset C^{\mathrm{Lip}} (X)$ 
such that $f_i \to f$ in $W^{1,2} (X)$ as $i \to \infty$. 
Then \eqref{eq:BE} yields $T_t f_i \in C^{\mathrm{Lip}} (X)$. 
Let $y \in X$ and 
$\gamma \: : \: [0,l] \to X$  
a unit speed minimal geodesic from $x$ to $y$. 
Since $\gl{T_t f_i}$ is an upper gradient of $T_t f_i$ 
(see \cite[Proposition~1.11]{Chee99} for instance), 
we have 
\begin{align} \nonumber
| T_t f_i (y) - T_t f_i (x) | 
& \le 
\int_0^l \gl{T_t f_i} ( \gamma (s) ) \,ds 
\\ \label{eq:BE-Lip}
& \le 
\e^{-Kt} \int_0^l T_t ( \gl{f_i}^2 )( \gamma (s) )^{1/2} \,ds , 
\end{align}
where the second inequality follows from 
\eqref{eq:BE}. 
Thanks to \eqref{eq:D} and the boundedness of $p_t$, 
$T_t f_i$ converges pointwisely to $T_t f$ as $i \to \infty$.  
Since $T_t ( \gl{f_i}^2 )(z) = T_t ( \gs{f_i}^2 ) (z)$ 
for $z \in X$, 
$T_t ( \gl{f_i}^2 )$ converges pointwisely 
to $T_t ( \gs{f}^2 )$ in a similar manner. 
Thus, by letting $i \to \infty$ in \eqref{eq:BE-Lip}, 
we obtain 
\begin{equation} \label{eq:BE-Lip2}
| T_t f (x) - T_t f (y) | 
  \le 
\e^{-Kt} 
\int_0^l T_t ( \gs{f}^2 ) ( \gamma (s) )^{1/2} \,ds . 
\end{equation}
The boundedness of $p_t$ yields that 
there is $C > 0$ satisfying 
\begin{equation} \label{eq:ultra}
T_t g  ( z )   
  \le 
C \| g \|_{L^1 (X,\vol)} 
\end{equation}
for all $g \in L^1 (X, \vol)$ and $z \in X$. 
Since $l = d( x, y )$, 
the estimate \eqref{eq:ultra} for $g = \gs{f}^2$ 
together with \eqref{eq:BE-Lip2} 
implies $T_t f \in C^{\mathrm{Lip}} (X)$. 
In order to show \eqref{eq:BE2}, 
choose a sequence $\{ y_i \}_{i \in \N}$ in $X$ so that 
$y_i \to x$ as $i \to \infty$ and 
\begin{equation*}
\lim_{i \to \infty} 
\frac{ | T_t f (x) - T_t f ( y_i ) | }{ d ( x, y_i ) } 
= \gl{T_t f} (x) .
\end{equation*}
By the continuity of $p_t$, $T_t ( \gs{f}^2 ) \in C (X)$ holds. 
Thus, applying \eqref{eq:BE-Lip2} for $y = y_i$, 
dividing both sides of it by $d( x, y_i )$ 
and letting $i \to \infty$ yield \eqref{eq:BE2}. 
\end{Proof}
As an easy but important consequence of Theorem~\ref{thm:BE}, 
we obtain the Lipschitz continuity of 
the heat kernel $p_t (x,y)$ 
as well as that of eigenfunctions. 
Recall that $-\lap$ has discrete spectrum 
consisting of nonnegative eigenvalues 
with finite multiplicity 
(\cite[Corollary~1.1]{Kuw_Mac_Shi}). 
\begin{theorem}[Lipschitz continuity] \label{thm:HK-Lip} 
\begin{enumerate}
\item
For $\nu \in \prob{X}$ and $t>0$, let $f_t = d T_t \nu / d \vol$.
Then $f_t \in C^{\mathrm{Lip}} (X)$ for $t > 0$. 
In particular, we have 
$p_t ( x, \cdot ) \in C^{\mathrm{Lip}} (X)$ and $T_t f \in C^{\mathrm{Lip}} (X)$ 
for all $x \in X$ and $f \in L^1 ( X, \vol )$. 
\item
Let $f$ be an $L^2$-eigenfunction of $\lap$ 
corresponding to an eigenvalue $- \lambda < 0$. 
Then $f \in C^{\mathrm{Lip}} (X)$. 
Moreover, 
$\gl{f} \le \e^{(\lambda - K ) t } \sqrt{\lambda \| T_t \|_{L^1 \to L^\infty}}
\| f \|_{L^2 (X , \vol)}$ holds for each $t > 0$. 
\end{enumerate}
\end{theorem}
\begin{Proof}
Since $f_t \in W^{1,2} (X)$, 
Theorem~\ref{thm:BE} yields that
$
f_t 
  = 
T_{t/2} f_{t/2} 
  \in 
C^{\mathrm{Lip}} (X)
$. 
For the second assertion, 
note that $f = \e^{\lambda t } T_t f$.  
Then the first assertion 
and $\mathcal{E} (f,f) = \lambda \| f \|_{L^2 (X, \vol)}^2$ 
yield the conclusion. 
\end{Proof}
\begin{remark} 
(i)
It has been known that the heat kernel $p_t (x,y)$ is
(locally) H\"older continuous of some fractional exponent,
that follows from the parabolic Harnack inequality 
shown in \cite{Kuw_Mac_Shi}. 

(ii)
The existence and the continuity of $p_t$ are used 
in the proof of Theorems~\ref{thm:BE}, \ref{thm:HK-Lip}
in an essential way (cf. Remark~\ref{rem:initial}). 

(iii)
To obtain a useful estimate of $\gl{T_t f}$ along our argument, 
we need a nice bound for $\| T_t \|_{L^1 \to L^\infty}$. 
For instance, 
the parabolic Harnack inequality, or the Nash inequality, 
implies 
\begin{equation}\label{eq:DHK}
p_t (x,y) \le \frac{C}{ \vol ( B_{\sqrt{t}} (x)) }
\end{equation} 
with some constant $C > 0$ being independent of $t,x,y$, for small $t$. 
It gives a bound for $\| T_t \|_{L^1 \to L^\infty}$. 
By a general argument, \eqref{eq:DHK} follows from 
the local Poincar\'e inequality and the volume doubling condition, 
both of which depend
only on the dimension $n$ of $X$ and a lower curvature bound
(see \cite{Kuw_Mac_Shi,Kuw_Shi01,Ranj_PI,Sturm_Harnack} for instance). 
However, we should be careful 
if we want to know whether $C$ in \eqref{eq:DHK} depends 
on the diameter and/or the volume of $X$. 
Indeed, estimates of type \eqref{eq:DHK} are mainly studied 
on non-compact state spaces and hence 
they did not seem to pay so much attentions 
on such a dependency in the literature. 
\end{remark}
In what follows, we consider two additional applications of 
the Bakry-\'Emery gradient estimate \eqref{eq:BE3} 
by employing the Lipschitz continuity of $T_t f$. 
The first one is the following inequality: 
\begin{theorem}[$\Gamma_2$-condition] \label{thm:Boch}
Let $f \in D ( \lap )$ with $\lap f \in W^{1,2} (X)$.  
Then, for $g \in D (\lap) \cap L^\infty ( X , \vol )$ 
with $g \ge 0$ and $\lap g \in L^\infty ( X, \vol )$, 
we have 
\begin{equation} \label{eq:Boch}
\frac12 
\int_X \lap g \gss{f}{f} d \vol 
- 
\int_X g \gss{\lap f}{f} d \vol  
  \ge 
K \int_X g \gss{f}{f} d \vol . 
\end{equation}  
\end{theorem}
\begin{remark} \label{rem:Boch}
(i) 
By virtue of the analyticity of $T_t$ 
and Theorem~\ref{thm:HK-Lip}, 
we have 
$T_t f \in D ( \lap^{m+1} )$ 
and $\lap^m T_t f \in C^{\mathrm{Lip}} (X)$ 
for any $f \in L^2 ( X , \vol )$, $t>0$ and $m \ge 0$. 
Thus there are fairly many $f$ and $g$'s 
satisfying the condition in Theorem~\ref{thm:Boch}. 

(ii) 
The inequality \eqref{eq:Boch} is nothing but a weak form of 
the $\Gamma_2$-condition 
\begin{equation} \label{eq:G2}
\Gamma_2 (f , f) 
: = 
\frac12
\big\{
  \lap ( \gss{f}{f} ) 
  - 2 \gss{f}{ \lap f} \big\}
\ge 
K \gss{f}{f} .
\end{equation}
This inequality is known to be 
equivalent to \eqref{eq:BE3} 
in an abstract framework 
(see \cite{Bak97,Led_geom-Markov} 
and references therein, for instance). 
However, the assumption involves 
the existence of a nice core 
$\mathcal{A} \subset D ( \lap )$ 
and it seems hopeless 
to verify it on Alexandrov spaces. 
When $X$ is a complete Riemannian manifold, 
the inequality \eqref{eq:G2} 
is equivalent to $\Ric \ge K$. 
Indeed, 
\begin{equation*}
\Gamma_2 (f,f) 
  =  
\Ric ( \nabla f , \nabla f ) 
  + 
| \mathop{\mathrm{Hess}} f |^2 
\end{equation*}
holds by the Bochner identity. 
\end{remark}
\begin{Proof}
We first show the claim 
for $f \in D ( \lap ) \cap C^{\mathrm{Lip}} ( X )$ 
with 
$\lap f \in D ( \lap ) \cap L^\infty (X , \vol)$.  
By \eqref{eq:BE3}, we obtain 
\begin{equation} \label{eq:Boch-1}
\int_X g \gss{ T_t f }{ T_t f } d \vol 
  \le 
\e^{-2Kt} 
\int_X g T_t ( \gss{f}{f} ) \,d \vol . 
\end{equation}
The derivation property $(\ref{eq:chain})$ yields that, 
for $t \ge 0$,    
\begin{align} \nonumber
\int_X g \gss{ T_t f }{ T_t f } d \vol 
& = 
\int_X \gss{ ( g T_t f ) }{ T_t f } d \vol 
  - 
\int_X T_t f \gss{g}{ T_t f } d \vol 
\\ \nonumber
& = 
  - 
\int_X g T_t f \lap T_t f \, d \vol  
  - 
\frac12 \int_X \gss{g}{ ( T_t f )^2 } d \vol 
\\ \label{eq:Gamma0}
& = 
  - 
\int_X g T_t f T_t \lap f \,d \vol 
  +  
\frac12 \int_X \lap g ( T_t f )^2 \,d \vol 
. 
\end{align}
Hence we obtain 
\begin{align} \nonumber 
\frac{d}{d t} 
& 
\left. 
  \int_X g \gss{ T_t f }{ T_t f } d \vol 
\right|_{t = 0}
\\ \nonumber
&  = 
- \int_X g ( \lap f )^2 d \vol 
- \int_X g f \lap^2 f d \vol 
+ \int_X \lap g f \lap f \, d \vol 
\\ \nonumber
& = 
\int_X \gss{ (g \lap f) }{ f } d \vol 
+ \int_X \gss{ (g f) }{ \lap f } d \vol 
- \int_X \gss{ g }{ (f \lap f) } d \vol 
\\ \label{eq:Boch-L}
& = 
2 \int_X g \gss{ \lap f }{ f } d \vol 
\end{align}
by using the derivation property again. 
Since $f \in C^{\mathrm{Lip}} (X)$, 
$\gss{f}{f} = \gl{f}^2 \in L^\infty (X ,\vol)$ holds.   
Hence we have 
\begin{align} \nonumber 
\frac{d}{d t}
& 
\left. 
  \abra{ 
    \e^{-2Kt} 
    \int_X g T_t ( \gss{f}{f} ) \,d \vol 
  }
\right|_{t=0} 
\\ \nonumber
& = 
\left. 
  \frac{d}{dt} 
  \abra{ 
    \e^{-2Kt} 
    \int_X (T_t g) \gss{f}{f} d \vol  
  }
\right|_{t=0} 
\\ \label{eq:Boch-R}
& = 
\int_X \lap g \gss{f}{f} d \vol
- 2 K \int_X g \gss{f}{f} d \vol . 
\end{align}
Since \eqref{eq:Boch-1} implies 
\begin{equation} \label{eq:Boch0}
\left. 
  \frac{d}{dt} 
  \int_X g \gss{ T_t f }{ T_t f } d \vol 
\right|_{t = 0}
  \le 
\left. 
  \frac{d}{dt}
  \abra{ 
    \e^{-2Kt} 
    \int_X g T_t ( \gss{f}{f} ) \,d \vol  
  }
\right|_{t=0} , 
\end{equation} 
we obtain \eqref{eq:Boch} 
by combining \eqref{eq:Boch0} 
with \eqref{eq:Boch-L} and \eqref{eq:Boch-R}. 

Next we consider the case that 
$f \in D (\lap)$ 
with $\lap f \in W^{1,2}(X)$. 
Then, by the above discussion, 
$T_\delta f$ and $g$ satisfy \eqref{eq:Boch} 
for $\delta > 0$ 
(cf.\ Remark~\ref{rem:Boch}(i)). 
Since 
$g, \lap g \in L^\infty ( X , \vol )$ 
and 
$\lap T_\delta f = T_\delta \lap f$, 
it suffices to show the claim that 
$\lim_{\delta \to 0} \gss{T_\delta h}{T_\delta h'} = \gss{h}{h'}$ 
weakly in $L^1 ( X, \vol )$ for $h, h' \in W^{1,2}(X)$. 
By polarization, we may assume $h = h'$. 
The spectral decomposition yields 
\begin{equation} \label{eq:app1} 
\lim_{\delta \to 0} 
\mathcal{E} ( T_\delta h - h , T_\delta h - h ) 
= 0
\end{equation}
(see \cite[Lemma~1.3.3]{FOT}, for instance). 
Let $\psi \in L^\infty (X, \vol)$. 
Then the Schwarz inequality yields 
\begin{align} \nonumber
&\Bigg| 
\int_X
\psi 
\gss{T_\delta h}{T_\delta h} d \vol
  - 
\int_X \psi \gss{h}{h}
d \vol 
\Bigg|
\\ \nonumber
& \le 
\abra{ 
  \abra{ 
    \int_X 
      \psi^2 \gs{T_\delta h}^2 
    d \vol
  }^{1/2} 
  + 
  \abra{ 
    \int_X 
      \psi^2 \gs{h}^2 
    d \vol
  }^{1/2} 
}
\mathcal{E} ( T_\delta h - h , T_\delta h - h )^{1/2}
\\ \label{eq:app2}
& \le 
2 \| \psi \|_{L^\infty}
\mathcal{E} ( h , h )^{1/2} 
\mathcal{E} ( T_\delta h - h , T_\delta h - h )^{1/2}.
\end{align}
Hence the desired claim follows 
from \eqref{eq:app1} and \eqref{eq:app2}.  
\end{Proof}
While we proved the implication 
from \eqref{eq:contraction2} with $p=2$ 
to \eqref{eq:BE3} and \eqref{eq:Boch}, 
these conditions are equivalent to each other 
on complete Riemannian manifolds (see \cite{vRS}). 
Such an equivalence still holds in our framework 
with a sharp constant, which can be different from $K$ in our hypothesis $\CD(K,\infty)$:
\begin{theorem}[Equivalence of ``Ricci curvature bound'' inequalities]\label{thm:Ricci}
Given $K_0 \in \R$, the following are equivalent. 
\begin{enumerate}
\item
\eqref{eq:contraction2} holds 
for $\nu_0 , \nu_1 \in \prob{X}$ and $t \ge 0$, 
with $p=2$ and $K = K_0$.  
\item
\eqref{eq:BE3} holds for $f \in W^{1,2} (X)$ and $t \ge 0$ 
with $K = K_0$. 
\item
\eqref{eq:Boch} holds with $K = K_0$ 
for 
$f \in D ( \lap )$ with $\lap f \in W^{1,2} (X)$ 
and $g \in D (\lap) \cap L^\infty ( X , \vol )$ 
with $g \ge 0$ and $\lap g \in L^\infty ( X, \vol )$.  
\end{enumerate}
\end{theorem}
\begin{Proof}
``(i) $\Rightarrow$ (ii) $\Rightarrow$ (iii)'' 
is already shown 
in Theorem~\ref{thm:BE} and Theorem~\ref{thm:Boch}. 

For ``(iii) $\Rightarrow$ (ii)'', it follows 
from a standard argument of 
the so-called $\Gamma_2$-calculus 
(see \cite{Bak97,Led_geom-Markov} for instance). 
For completeness, we give a sketch of the proof. 
Take $g_0 \in C (X)$ 
with $g_0 \ge 0$ arbitrary 
and let $g = T_\delta g_0$. 
Then $g \in D (\lap) \cap L^\infty ( X , \vol )$ 
with $g \ge 0$ and $\lap g \in L^\infty ( X, \vol )$. 
Let us define $\Psi : [0,t] \to \R$ by 
\[
\Psi (s) := \int_X g T_s ( \gs{T_{t-s} f}^2 ) \,d \vol .
\]
By a similar calculation as in \eqref{eq:Gamma0}, 
we can easily prove that 
$\Psi$ is continuous on $[0,t)$ and   
$C^1$ on $( 0 , t )$. 
A similar argument as in \eqref{eq:app2} yields 
that $\Psi$ is continuous at $t$. 
Here we use 
the ultracontractivity $\| T_t \|_{L^1 \to L^\infty} < \infty$. 
A similar calculation 
as in the proof of Theorem~\ref{thm:BE} 
together with the assumption in (iii) 
leads to the inequality 
$\Psi' (s) \ge 2 K_0 \Psi(s)$. 
Hence \eqref{eq:BE3} follows by integrating it. 

For ``(ii) $\Rightarrow$ (i)'', we claim 
that \eqref{eq:BE3} implies 
\eqref{eq:BE2} for \emph{every} $x \in X$. 
Indeed, by using a bi-Lipschitz chart, we can bring the problem 
locally on an open set in a Euclidean space. 
Then, by applying \cite[Lemma~3.2.1]{AT} and 
by using the continuity of $T_t ( \gs{f}^2 )$, 
the claim follows. 
Then we can apply \cite[Theorem~2.2]{K9} to conclude (i) 
from \eqref{eq:BE2}. 
\end{Proof}

\begin{remark}
By the same argument 
as in ``(ii) $\Rightarrow$ (i)'' of the last proof, 
we can give a proof of $T_t f \in C^{\mathrm{Lip}} (X)$ 
for $f \in W^{1,2} (X)$ 
under the condition \eqref{eq:BE3}. 
In other word, 
a priori regularity 
$T_t f \in C^{\mathrm{Lip}} (X)$ by Theorem~\ref{thm:BE} 
is not used in the last proof. 
\end{remark}
As the second application of \eqref{eq:BE3}, 
we mention that \eqref{eq:BE3} 
together with Theorem~\ref{thm:HK-Lip} 
implies some functional inequalities 
by means of \cite[Theorem~1.3]{Kaw-Miyo}. 
Since $T_t$ is Markovian, 
we can restrict $T_t$ 
to a contraction on $L^\infty ( X, \vol )$.
Then we can further extend $T_t$ to a contraction 
on $L^p (X, \vol)$ for any $1 \le p \le \infty$ 
by interpolation and the symmetry of $T_t$. 
Let us denote the infinitesimal generator 
of $T_t$ in $L^p ( X , \vol )$ 
by $\lap_p$. 
Let us define 
$R_\alpha^{(q)} f : = \gs{ ( ( \alpha - \lap_p )^{-q/2} f ) }$. 
\begin{corollary} 
Let $2 \le p < \infty$, $q > 1$ and 
$\alpha > \max\{ ( - K ) , 0 \}$. 
Then we have the following: 
\begin{enumerate}
\item
There exists a constant $C_R > 0$ which depends 
only on $p,q$ and $\max\{ ( \alpha + K ) , \alpha \}$ 
such that 
\begin{equation*}
\| 
  R_\alpha^{(q)} f 
\|_{L^p} 
  \le 
C_R 
\| f \|_{L^p}
\end{equation*}
for $f \in L^p ( X, \vol )$. 
\item 
Suppose $q < 2$. 
Then there exists $C_{p,q} > 0$ 
such that  
\begin{equation*}
\left\| \gs{T_t f} \right\|_{L^p} 
  \le 
C_{p,q} 
\| R_\alpha^{(q)} \|_{L^p \to L^p} 
\abra{ 
  \alpha^{q/2} + t^{-q/2} 
}
\| f \|_{L^p} 
\end{equation*}
for $t > 0$ and $f \in L^p ( X , \vol )$. 
\end{enumerate}
\end{corollary} 
\begin{Proof}
It is sufficient to verify 
the assumption of \cite[Theorem~1.3]{Kaw-Miyo}. 
Set 
$
\mathcal{A} 
  : = 
C^{\mathrm{Lip}} (X) \cap D ( \lap_2 )
$. 
Since we already know \eqref{eq:BE3}, 
we only need to show the following claim: 
$\mathcal{A}$ is dense 
in $W^{1,2} (X)$ 
and 
$f^2 \in D ( \lap_1 )$ holds 
for any $f \in \mathcal{A}$. 
By Theorem~\ref{thm:HK-Lip}, 
$T_t f \in \mathcal{A}$ holds. 
Thus $\mathcal{A}$ is dense in $W^{1,2} (X)$ 
since $\mathcal{E} ( T_t f - f , T_t f - f ) \to 0$ 
as $t \to 0$ (cf.\ \eqref{eq:app1}). 
By \cite[Proposition~I.2.4.3]{BH}, 
it is enough to prove $f^2 \in D ( \lap_2 )$ 
for $f \in \mathcal{A}$. 
Take $g \in W^{1,2} (X)$ arbitrary. 
Recall that $\gs{f}^2 = \gl{f}^2 \in L^\infty (X, \vol)$. 
The derivation property yields 
\begin{align*}
| \mathcal{E} (f^2 , g) |
&  = 
2 
\left| 
  \int_X f \gss{f}{g} d \vol 
\right| 
\\
& = 
2 
\left| 
  \int_X 
    ( \gss{f}{(gf)} - g \gss{f}{f})
  d \vol 
\right|
\\
& = 
2 
\left| 
  \int_X 
    ( f \lap f - \gss{f}{f} ) g 
  \, d \vol 
\right|
\\
& \le 
2 
\abra{ 
  \| f \|_{L^\infty}
  \| \lap f \|_{L^2} 
    + 
  \| \gs{f}^2 \|_{L^2}
} \| g \|_{L^2}. 
\end{align*}
This estimate means $f^2 \in D ( \lap_2 )$ 
and hence the proof is completed.  
\end{Proof}

Finally, we observe that all of our results are generalized to
the heat equation with drift, in other words, the Fokker-Planck equation. 
Given a potential function $V \in C^{\mathrm{Lip}} (X)$,
we modify the Dirichlet energy and the relative entropy into
\begin{align*}
\mathcal{E}^V ( f, g ) 
& = 
\int_X \gss{f}{g} \e^{-V} d\vol \qquad f,g \in W^{1,2}(X),
\\
\Ent^V ( \mu ) 
& = 
\Ent ( \mu ) + \int_X V \, d \mu \qquad \mu \in \prob{X}.
\end{align*}
We regard $\mathcal{E}^V$ 
as a bilinear form on $L^2 ( X, \e^{-V} \vol )$. 
Observe that 
$\Ent^V$ is nothing but the relative entropy with respect to 
$\e^{-V} \vol$. 
Note that the semigroup $T_t^V$ on $L^2 ( X, \e^{-V} \vol )$ 
associated with $\mathcal{E}^V$ solves 
the following diffusion equation 
\begin{equation}\label{eq:FP}
\frac{d}{dt} u = \lap u -  \gss{V}{u}. 
\end{equation}
Since $\e^{-V}$ is bounded and away from 0, 
$\e^{-V} \vol$ is equivalent to $\vol$. 
Hence $(X, d , \e^{-V} \vol)$ satisfies the volume doubling condition
as well as  the Poincar\'e inequality for upper gradients. 
Moreover, as $\mathcal{E}$ and $\mathcal{E}^V$ are equivalent, 
the Poincar\'e inequality for $\mathcal{E}^V$ is also valid. 
Therefore a continuous density $p_t^V$ for $T_t^V$ exists. 
Under the assumption that $\Ent^V$ is $K$-geodesically convex, 
we can apply the general theory of the gradient flow on $( \prob{X} , W_2 )$
to obtain the gradient flow $\mu_t$ of $\Ent^V$. 
Furthermore, every argument in Sections~\ref{sc:main}, \ref{sc:appl}
works verbatim and gives similar results for the equation $(\ref{eq:FP})$. 
Note that, under $\CD(K,\infty)$ for $(X,d,\vol)$, the $K'$-convexity of $V$ 
implies the $(K+K')$-convexity of $\Ent^V$. 
\end{section}

\bibliographystyle{amsplain}

\begin{thebibliography}{99}

\bibitem{AGS}
L.~Ambrosio, N.~Gigli, and G.~Savar\'e,
\emph{Gradient flows in metric spaces and in the space of probability measures}.
Second edition, Birkh\"auser Verlag, Basel, 2008.

\bibitem{AT}
L.~Ambrosio and P.~Tilli, 
\emph{Topics on analysis in metric spaces}. 
Oxford University Press, Oxford, 2004. 

\bibitem{Bak97}
D.~Bakry, \emph{On {S}obolev and logarithmic {S}obolev inequalities for Markov
  semigroups}, New trends in stochastic analysis (Charingworth, 1994), World
  Sci.\ Publ.\ River Edge, NJ, 1997, pp.~43--75.

\bibitem{BEHM_HJ}
Z.~M.~Balogh, A.~Engoulatov, L.~Hunziker, and O.~E.~Maasalo, \emph{Functional
  inequalities and {H}amilton-{J}acobi equations in geodesic spaces},
to appear in Potential Anal.

\bibitem{Bert}
J.~Bertrand, \emph{Existence and uniqueness of optimal maps on Alexandrov spaces},
Adv.\ Math.\ {\bf 219} (2008), 838--851.

\bibitem{BH}
N.~Bouleau and F.~Hirsch, \emph{Dirichlet forms and analysis on {W}iener
  space}, de Gruyter Studies in Mathematics, {\bf 14}, Walter de Gruyter \& Co.,
  Berlin, 1991.

\bibitem{BBI}
D.~Burago, Yu. Burago, and S.~Ivanov, \emph{A course in metric geometry},
  Graduate studies in mathematics, {\bf 33}, American mathematical society,
  Providence, RI, 2001.

\bibitem{BGP}
Yu.~Burago, M.~Gromov, and G.~Perel'man,
\emph{A.\ D.\ Alexandrov spaces with curvatures bounded below} (Russian),
Uspekhi Mat.\ Nauk {\bf 47} (1992), 3--51, 222;
English translation: Russian Math.\ Surveys {\bf 47} (1992), 1--58.

\bibitem{Chee99}
J.~Cheeger, \emph{Differentiability of {L}ipschitz functions on metric measure
  spaces}, Geom.\ Funct.\ Anal.\ \textbf{9} (1999), no.~3, 428--517.

\bibitem{FSS}
S.~Fang, J.~Shao, and K.-T.~Sturm,
\emph{Wasserstein space over the Wiener space},
Probab.\ Theory Related Fields {\bf 146} (2010), no.~3-4, 535--565.

\bibitem{FJ}
A.~Figalli and N.~Juillet,
\emph{Absolute continuity of Wasserstein geodesics in the Heisenberg group},
J.\ Funct.\ Anal.\ {\bf 255} (2008), no.~1, 133--141.

\bibitem{FOT}
M.~Fukushima, Y.~Oshima, and M.~Takeda, \emph{Dirichlet forms and symmetric
  {M}arkov processes}, de Gruyter Studies in Mathematics, {\bf 19}, Walter de Gruyter
  \& Co., Berlin, 1994.

\bibitem{heatmm}
N.~Gigli,
\emph{On the heat flow on metric measure spaces: existence, uniqueness and stability},
Calc.\ Var.\ Partial Differential Equations {\bf 39} (2010), 101--120.

\bibitem{Gigli-Ohta10}
N.~Gigli and S.~Ohta,
\emph{First variation formula in {W}asserstein spaces
 over compact {A}lexandrov spaces},
to appear in Canad.\ Math.\ Bull.

\bibitem{JKO:FP:98}
R.~Jordan, D.~Kinderlehrer, and F.~Otto,
\emph{The variational formulation of the {F}okker-{P}lanck equation}.
SIAM J.\ Math.\ Anal.\ {\bf 29} (1998), 1--17.

\bibitem{Jui}
N.~Juillet, \emph{Diffusion by optimal transport in Heisenberg groups},
Preprint. Available at {\sf http://www-irma.u-strasbg.fr/\textasciitilde juillet/}

\bibitem{Kaw-Miyo}
H.~Kawabi and T.~Miyokawa, \emph{The {L}ittlewood-{P}aley-{S}tein inequality
  for diffusion processes on general metric spaces}, J.\ Math.\ Sci.\ Univ.\ Tokyo
  \textbf{14} (2007), no.~1, 1--30.

\bibitem{Kor-Sch}
N.~J.~Korevaar and R.~M.~Schoen,
\emph{Sobolev spaces and harmonic maps for metric space targets},
Comm.\ Anal.\ Geom.\ {\bf 1} (1993), 561--659.

\bibitem{K9}
K.~Kuwada, \emph{Duality on gradient estimates and {W}asserstein controls}, J.\
  Funct.\ Anal.\ \textbf{258} (2010), no.~11, 3758--3774.

\bibitem{Kuw_Mac_Shi}
K.~Kuwae, Y.~Machigashira, and T.~Shioya, \emph{Sobolev spaces, {L}aplacian and
  heat kernel on {A}lexandrov spaces}, Math.\ Z.\ \textbf{238} (2001), no.~2,
  269--316.

\bibitem{Kuw_Shi01}
K.~Kuwae and T.~Shioya, \emph{On generalized measure contraction property and
  energy functionals over {L}ipschitz maps}, Potential Anal.\ \textbf{15}
  (2001), no.~1, 105--121.

\bibitem{Kuw_Shi03}
\bysame, \emph{Sobolev spaces and {D}irichlet spaces over maps between metric
  spaces}, J.\ Reine Angew.\ Math.\ \textbf{555} (2003), 39--75.

\bibitem{Led_geom-Markov}
M.~Ledoux, \emph{The geometry of {M}arkov diffusion generators}, Ann.\ Fac.\ Sci.\
  Toulouse Math.\ (6) \textbf{9} (2000), no.~2, 305--366.

\bibitem{LV-CD2}
J.~Lott and C.~Villani,
\emph{Weak curvature conditions and functional inequalities},
J.\ Funct.\ Anal.\ {\bf 245} (2007), 311--333.

\bibitem{HJ}
\bysame,
\emph{Hamilton-Jacobi semigroup in length spaces and applications},
J.\ Math.\ Pures Appl.\ {\bf 88}, (2007), 219--229.

\bibitem{LV-CD}
\bysame, \emph{Ricci curvature for metric-measure spaces via optimal transport},
Ann.\ of Math.\ {\bf 169} (2009), 903--991.

\bibitem{Ohta_grad-Alex}
S.~Ohta, \emph{Gradient flows on {W}asserstein spaces over compact {A}lexandrov
  spaces}, Amer. J. Math \textbf{131} (2009), no.~2, 475--516.

\bibitem{heatfinsl}
S.~Ohta and K.-T.~Sturm, \emph{Heat flow on Finsler manifolds},
Comm.\ Pure.\ Appl.\ Math.\ {\bf 62} (2009), 1386--1433.

\bibitem{OS2}
\bysame,
\emph{Non-contraction of heat flow on Minkowski spaces},
to appear in Arch.\ Ration.\ Mech.\ Anal.

\bibitem{Otsu-Shioya_JDG94}
Y.~Otsu and T.~Shioya, \emph{The {R}iemannian structure of {A}lexandrov
  spaces}, J. Differential Geom. \textbf{39} (1994), no.~3, 629--658.

\bibitem{otto:geom:01}
F.~Otto,
\emph{The geometry of dissipative evolution equations: the porous medium equation},
Comm.\ Partial Differential Equations {\bf 26} (2001), no.\ 1--2, 101--174.

\bibitem{Perel}
G.~Perelman,
\emph{DC structure on Alexandrov space with curvature bounded below},
Available at {\sf http://www.math.psu.edu/petrunin/}

\bibitem{PetERA}
A.~Petrunin,
\emph{Harmonic functions on Alexandrov spaces and their applications},
Electron.\ Res.\ Announc.\ Amer.\ Math.\ Soc.\ {\bf 9} (2003), 135--141. 

\bibitem{PetCD}
\bysame, \emph{Alexandrov meets Lott--Villani--Sturm},
M\"unster J.\ Math.\ {\bf 4} (2011), 53--64.

\bibitem{QZZ}
Z.~Qian, H.-C.~Zhang and X.-P.~Zhu,
\emph{Sharp spectral gap and Li-Yau's estimate on Alexandrov spaces},
Preprint (2011). Available at {\sf arXiv:1102.4159}

\bibitem{Ranj_PI}
A.~Ranjbar-Motlagh, \emph{Poincar{\'e} inequality for abstract spaces},
Bull.\ Austral.\ Math.\ Soc.\ \textbf{71} (2005), no.~2, 193--204.

\bibitem{vRS}
M.-K.~von Renesse and K.-T.~Sturm,
\emph{Transport inequalities, gradient estimates, entropy and Ricci curvature},
Comm.\ Pure Appl.\ Math.\ {\bf 58} (2005), 923--940.

\bibitem{Sav}
G.~Savar\'e,
\emph{Gradient flows and diffusion semigroups in metric spaces
under lower curvature bounds},
C.\ R.\ Math.\ Acad.\ Sci.\ Paris {\bf 345} (2007), 151--154.

\bibitem{Sturm_Harnack}
K.-T.~Sturm, \emph{Analysis on local {D}irichlet spaces.{III}. {T}he parabolic
  {H}arnack inequality}, J.\ Math.\ Pures Appl.\ (9) \textbf{75} (1996), no.~3,
  273--297.

\bibitem{SturmI}
\bysame, \emph{On the geometry of metric measure spaces.~I},
Acta Math.\ {\bf 196} (2006), 65--131.

\bibitem{SturmII}
\bysame, \emph{On the geometry of metric measure spaces.~II},
Acta Math.\ {\bf 196} (2006), 133--177.

\bibitem{Villani08}
C.~Villani, \emph{Optimal Transport, old and new}, Grundlehren der mathematischen Wissenschaften collection, 2009.

\bibitem{ZZ1}
H.-C.~Zhang and X.-P.~Zhu
\emph{Ricci curvature on Alexandrov spaces and rigidity theorems},
Comm.\ Anal.\ Geom.\ {\bf 18} (2010), no.~3, 503--553.

\bibitem{ZZ2}
\bysame, \emph{On a new definition of Ricci curvature on Alexandrov spaces},
Acta.\ Math.\ Sci.\ Ser.\ B Engl.\ Ed.\ {\bf 30} (2010), no.\ 6, 1949--1974.

\bibitem{ZZ3}
\bysame, \emph{Yau's gradient estimate on Alexandrov spaces}, 
Preprint (2010). Available at {\sf arXiv:1012.4233} 
\end{thebibliography}

\providecommand{\bysame}{\leavevmode\hbox to3em{\hrulefill}\thinspace}
\providecommand{\MR}{\relax\ifhmode\unskip\space\fi MR }
\providecommand{\MRhref}[2]{%
  \href{http://www.ams.org/mathscinet-getitem?mr=#1}{#2}
}
\providecommand{\href}[2]{#2}

\end{document}